\documentclass[12pt]{amsart}
\usepackage{amssymb}

\newtheorem{thm}{Theorem}[section]

\newtheorem{lem}[thm]{Lemma}
\newtheorem{prop}[thm]{Proposition}

\theoremstyle{remark}

\newtheorem{rem}[thm]{Remark}
\newtheorem{dfn}[thm]{Definition}
\newtheorem{thma}[thm]{Theorem}
\newtheorem{propa}[thm]{Proposition}
\newtheorem{lema}[thm]{Lemma}

\def\Aut{{\rm Aut}}
\def\Syl{{\rm Syl}}
\def\id{{\rm id}}

\def\Out{{\rm Out}}
\def\Hom{{\rm Hom}}

\def\N#1#2{N_{#1}(#2)}
\def\C#1#2{C_{#1}(#2)}

\def\Cc{{\mathcal C}}

\def\Ff{{\mathcal F}}
\def\Ee{{\mathcal E}}
\def\Gg{{\mathcal G}}
\def\Uu{{\mathcal U}}

\begin{document}

\title[A characteristic subgroup for fusion systems]{A characteristic subgroup for fusion systems}

\author{Silvia Onofrei}
\author{Radu Stancu}
\address{Department of Mathematics, The Ohio State University, 100 Mathematics Tower, 231 West $18^{\rm th}$ Avenue, Columbus, Ohio 43210, USA, Email: onofrei@math.ohio-state.edu}
\address{Department of Mathematical Sciences, University of Copenhagen
Universitetsparken 5, DK-2100 Copenhagen {\O}, Denmark, Email: stancu@math.ku.dk}
\date{\today}

\begin{abstract}As a counterpart for the prime $2$ to Glauberman's $ZJ$-theorem, Stellmacher proves that any nontrivial $2$-group $S$ has a nontrivial characteristic subgroup $W(S)$ with the following property. For any finite $\Sigma_4$-free group $G$, with $S$ a Sylow $2$-subgroup of $G$ and with $O_2(G)$ self-centralizing, the subgroup $W(S)$ is normal in $G$. We generalize Stellmacher's result to fusion systems. A similar construction of $W(S)$ can be done for odd primes and gives rise to a Glauberman functor.
\end{abstract}

\maketitle

%%%%%%%%%%%%%%%%%%%%%%%%
\section{Introduction}%%
%%%%%%%%%%%%%%%%%%%%%%%%

A fundamental result in the theory of finite groups is Glauberman's $ZJ$-theorem \cite{gl68}. For $p$ be an odd prime, $G$ a finite group and $S$ a Sylow $p$-subgroup of $G$, the $ZJ$-theorem asserts that the center of the Thompson group $Z(J(S))$ is normal in $G$ whenever $G$ is $Qd(p)$-free and $\C G{O_p(G)} \leq O_p(G)$. Recall that for a finite $p$-group $Q$, the Thompson subgroup $J(Q)$ is the subgroup generated by the abelian subgroups of $Q$ of largest order and that the group $Qd(p)= (\mathbb{Z}_p \times \mathbb{Z}_p):SL(2,p)$ is the extension of the $2$-dimensional vector space over $\mathbb{F}_p$ (the field with $p$ elements) by $SL(2,p)$ with its natural action on this vector space. A group $G$ is $H$-free if no section of $G$ is isomorphic to $H$; see also Section~$3$.

More recently, a proof of the $ZJ$-theorem, in the context of fusion systems, was given by Kessar and Linckelman \cite{kl05}. The authors introduce the notion of $Qd(p)$-free fusion system and prove that if $\Ff$ is a $Qd(p)$-free fusion system on a finite $p$-group $S$, with $p$ an odd prime, then $\Ff$ is controlled by $W(S)$, for any Glauberman functor $W$. The related notions of characteristic $p$-functor and Glauberman functor were initially defined in \cite[Definition 1.3]{klr}; they are given below in Definition \ref{charfunct}.

For $p=2$ the $ZJ$-theorem does not hold anymore. In \cite[Question 16.1]{gl71}, Glauberman asks whether there exists a subgroup which is characteristic in a Sylow $2$-subgroup $S$ of a $\Sigma_4$-free group $G$, with the property $\C G{O_p(G)} \leq O_p(G)$. Here $\Sigma_4$ denotes the symmetric group on four letters.

The answer to Glauberman's question was given by Stellmacher \cite{st96}, who also obtained a different proof of the $ZJ$-theorem \cite{st92,st92b}. Stellmacher's idea is to approximate such a subgroup via subgroups of $Z(J(S))$; see \cite[Section 9.4]{ks04} for an overview of this approach. The main theorem in \cite{st96} (see also 6.4 in the Appendix) can be phrased as follows:

\noindent{\it Theorem (Stellmacher): Let $S$ be a finite nontrivial $2$-group. Then there exists a nontrivial characteristic subgroup $W(S)$ of $S$ which is normal in $G$, for every finite $\Sigma_4$-free group $G$ with $S$ a Sylow $2$-subgroup and $\C G{O_2(G)} \leq O_2(G)$.}

Remark that the condition (III) in \cite{st96} is not necessary. A proof of this fact uses Lemmas 6.5 and 6.6 and Remark 6.7 in the Appendix.

In this paper we generalize Stellmacher's approach to fusion systems. Our main result is a proof of Stellmacher's version of the $ZJ$-theorem in the context of fusion systems:

\begin{thm}\label{mainS4}Let $S$ be a finite $2$-group and let $\Ff$ be a $\Sigma_4$-free fusion system over $S$. Then there exists a nontrivial characteristic subgroup $W(S)$ of $S$ with the property that $\Ff=\N \Ff {W(S)}$.
\end{thm}

Since Stellmacher's construction of $W(S)$ gives rise to a Glauberman functor (see Section $4$ for details) we can combine Theorem $B$ in \cite{kl05} with Theorem \ref{mainS4} in our paper to obtain the more general result which is independent of the nature of the prime $p$:

\begin{thm}\label{Qdpfree}Let $S$ be a finite $p$-group and let $\Ff$ be a $Qd(p)$-free fusion system over $S$. Then there exists a nontrivial characteristic subgroup $W(S)$ of $S$ with the property that $\Ff=\N \Ff {W(S)}$.
\end{thm}

Using the same construction for $W(S)$ as in the above theorem, the normal complement theorem due to Thompson \cite[9.4.7]{ks04} can be phrased as:

{\it Theorem (Thompson): Let $G$ be a finite group, $p$ an odd prime and $S$ a Sylow $p$-subgroup of $G$. Then $G$ has a normal $p$-complement provided $\N G{W(S)}$ has such a complement}.

Our third result generalizes Thompson's theorem to the class of fusion systems. This result is similar to Theorem $A$ in \cite{kl05}, except that we replace the group $Z(J)$ with the group $W(S)$ defined in Section~$4$:

\begin{thm}\label{pcomp}Let $\Ff$ be a fusion system over a finite $p$-group $S$, with $p$ an odd prime. Then $\Ff = \Ff_S(S)$ if and only if $\N \Ff {W(S)} = \Ff_S(S)$.
\end{thm}

The paper is organized as follows. Section $2$ contains background material on fusion systems. In Section $3$, the notions of $H$-free fusion system, characteristic $p$-functor and Glauberman functor are defined; further properties of fusion systems are discussed. The characteristic subgroup $W(S)$ is constructed, via two different methods, in Section $4$. The proofs of the theorems are given in Section $5$. In the Appendix a few related results from group theory are included.

%%%%%%%%%%%%%%%%%%%%%%%%%%%%%%%%%%%%%%%%
\section{Background on Fusion Systems}%%
%%%%%%%%%%%%%%%%%%%%%%%%%%%%%%%%%%%%%%%%

Fusion systems were introduced by Puig in 1990 \cite{puig:notes,puig:frobenius} in an effort to axiomatize the $p$-local structure of a finite group and of a block of a group algebra - the work was published only recently \cite{puig:fusioncategories} but was known to the community long before. In 2000 Broto, Levi and Oliver \cite{blo} 
enriched this axiomatic approach by introducing the centric linking system that is essentially linked to the associated $p$-completed classifying space. The third author used this axiomatic frame to reformulate and solve the Martino-Priddy conjecture \cite{oliver,oliver2}. Broto, Levi and Oliver give a different definition of the fusion systems which they proved to be equivalent to Puig's definition. In this paper we use a simplified definition which we find more elegant, equivalent to the above ones \cite{kessarstancu}.

We start with a more general definition, following \cite{linckelmann}.

A {\it category $\Ff$ on a finite $p$-group $S$} is a category whose objects are
the subgroups of $S$ and whose set of morphisms
between the subgroups $Q$ and $R$ of $S$, is a set $\Hom_\Ff(Q,R)$ of
injective group homomorphisms from $Q$ to $R$, with the following properties:
\begin{enumerate}
\item[(1)] if $Q\le R$ then the inclusion of $Q$ in $R$ is a morphism in $\Hom_\Ff(Q,R)$;
\item[(2)] for any $\varphi\in\Hom_\Ff(Q,R)$ the induced isomorphism $Q\simeq\varphi(Q)$
and its inverse are morphisms in $\Ff$;
\item[(3)] the composition of morphisms in $\Ff$ is the
usual composition of group homomorphisms.
\end{enumerate}

Let $\Ff_1$ be a category on $S_1$ and $\Ff_2$ a category on $S_2$.
A {\it morphism} between $\Ff_1$ and $\Ff_2$ is a pair $(\alpha, \Theta)$ with $\alpha \in \Aut(S)$ and
$\Theta: \Ff_1 \rightarrow \Ff_2$ a covariant functor, such that:
\begin{itemize}
\item[(i)]for any subgroup $Q$ of $S$, $\alpha (Q) = \Theta(Q)$;
\item[(ii)]for any morphism $\varphi$ in $\Ff_1$, $\Theta(\varphi)\circ \alpha = \alpha \circ \varphi$.
\end{itemize}

In the following we give a series of useful definitions in a category
$\Ff$ on $S$. If there exists an isomorphism $\varphi\in\Hom_\Ff(Q,R)$ we say that
$Q$ and $R$ are {\it $\Ff$-conjugate}.

We say that a subgroup $Q$ of $S$ is
\begin{list}{\upshape\bfseries}
{\setlength{\leftmargin}{.8cm}
\setlength{\rightmargin}{0cm}
\setlength{\labelwidth}{0.8cm}
\setlength{\labelsep}{0.2cm}
\setlength{\parsep}{0.5ex plus 0.1ex minus 0.1ex}
\setlength{\itemsep}{0.2ex plus 0.1ex minus 0.1ex}}
\item[(i)] {\it fully $\Ff$-centralized} if $|\C SQ|\ge|\C S{Q'}|$ for all $Q'\le S$
which are $\Ff$-conjugate to $Q$.
\item[(ii)] {\it fully $\Ff$-normalized} if $|\N SQ|\ge|\N S{Q'}|$ for all $Q'\le S$
which are $\Ff$-conjugate to $Q$.
\item[(iii)] {\it $\Ff$-centric} if $\C S{\varphi(Q)}\subseteq \varphi(Q)$, for all $\varphi\in\Hom_\Ff(Q,S)$.
\item[(iv)] {\it $\Ff$-radical} if $O_p(\Out_\Ff(Q))=1$.
\item[(v)] {\it $\Ff$-essential} if $Q$ is $\Ff$-centric and $\Out_\Ff(Q)$
has a strongly $p$-embedded proper subgroup $M$
(that is $M$ contains a Sylow $p$-subgroup $P$ of $\Out_\Ff(Q)$ such that
$P\ne\,^\varphi P$ and $^\varphi P\cap P=\{1\}$ for every $\varphi\in\,\Out_\Ff(Q)\setminus M$ ).
\end{list}

For $Q,\,R \leq S$ we denote $\Hom_S(Q,R):=\{u\in S\,|\,^uQ\le R\}/\C SQ$
and $\Aut_S(Q):=\Hom_S(Q,Q)$. Other useful notations are $\Aut_\Ff(Q):=\Hom_\Ff(Q,Q)$
and $\Out_\Ff(Q):=\Aut_\Ff(Q)/\Aut_Q(Q)$.

We are now ready to give the definition of a fusion system.

A {\it fusion system} on a finite $p$-group $S$ is a category $\Ff$ on $S$ satisfying the following
properties:
\begin{enumerate}
\item[FS1.] $\Hom_S(Q,R)\subseteq\Hom_\Ff(Q,R) \text{ for all } Q,R\le S\,.$

\item[FS2.] $\Aut_S(S)$ is a Sylow $p$-subgroup of $\Aut_\Ff(S)$.

\item[FS3.] Every $\varphi:Q\to S$ such that $\varphi(Q)$ is fully $\Ff$-normalized extends to a morphism $\widehat\varphi:N_\varphi\to S$ where
$$N_\varphi=\{x\in\N SQ\,|\,\exists y\in\N S{\varphi(Q)},\,\varphi(\,^xu)=\,^y\varphi(u),\,\forall u\in Q\}\,.$$
\end{enumerate}

Remark that $N_\varphi$ is the largest subgroup of $\N SQ$ such that
$^\varphi(N_\varphi/\C SQ)\le\Aut_S(\varphi(Q))$. Thus we always have $Q\C SQ\le N_\varphi \le N_S(Q)$.

If $\Ff$ is a fusion system on $S$ and $Q\le S$ we have the following
equivalent characterization of being fully $\Ff$-normalized.

\begin{prop}[\cite{markus}, Proposition 1.6]\label{prpsylow}
A subgroup $Q$ of $S$ is fully $\Ff$-normalized if and only if $Q$ is fully
$\Ff$-centralized and $\Aut_S(Q)$ is a Sylow $p$-subgroup of $\Aut_\Ff(Q)$.
\end{prop}

In the following Lemma we recollect two useful properties involving fully $\Ff$-normalized subgroups; see also \cite[Lemmas 2.2, 2.3]{kl05}. For completeness we include the proofs.

\begin{lem}\label{fullynormalized} Let $\Ff$ be a fusion system on a finite $p$-group $S$ and $Q$ a subgroup of $S$.
\begin{itemize}
\item[a)] There is a morphism $\varphi \in \Hom_\Ff(N_S(Q),S)$ such that $\varphi(Q)$ is fully $\Ff$-normalized.
\item[b)] If $Q$ is fully $\Ff$-normalized, then $\varphi(Q)$ is fully normalized, for any morphism $\varphi \in \Hom_\Ff(N_S(Q),S)$.
\end{itemize}
\end{lem}

\begin{proof}
a) Let $\psi:Q \rightarrow S$ be a morphism with $\psi(Q)$ fully $\Ff$-normalized. By Proposition $2.1$, $\Aut_S(\psi(Q))$ is a Sylow $p$-subgroup of $\Aut_\Ff (\psi(Q))$. Since $\psi \circ \Aut_S(Q) \circ \psi^{-1}$ is a $p$-subgroup of $\Aut_\Ff(\psi(Q))$ it follows that there exists a morphism $\tau \in \Aut_\Ff(\psi(Q))$ with $\tau \psi \circ \Aut_S(Q) \circ \psi^{-1} \tau^{-1} \le \Aut_S(\psi(Q))$. Set $\alpha = \tau \psi$ and observe that $\alpha(Q)$ is fully $\Ff$-normalized. By the extension axiom FS3, $\alpha$ extends to a morphism $\varphi:N_\alpha \rightarrow S$. But since $\alpha \circ \Aut_S(Q) \circ \alpha^{-1} \le \Aut_S(\psi(Q))$ it follows that $N_{\alpha} = N_S(Q)$. Henceforth there exists a morphism $\varphi \in \Hom _\Ff (\N S Q, S)$ such that $\varphi(Q)$ is fully $\Ff$-normalized.\\
b) Since $Q$ is fully $\Ff$-normalized and since $\varphi$ is a morphism in $\Ff$, hence injective, we have: $\varphi(N_S(Q))=N_S(\varphi(Q))$.
\end{proof}

Puig \cite{puig:notes} gave analogous notions for the normalizer and the centralizer in a finite group:

The {\it normalizer of $Q$ in $\Ff$} is the category $\N\Ff Q$ on $\N SQ$
having as morphisms those morphisms $\varphi\in\Hom_{\Ff}(R,T)$, for $R$ and $T$ subgroups of $\N SQ$, satisfying that there exists a morphism $\widehat{\varphi}\in\Hom_\Ff(QR,QT)$ such that $\widehat{\varphi}|_Q\in\Aut_\Ff(Q)$ and $\widehat{\varphi}|_R=\varphi$. If $Q \le S$ has the property that $\Ff=\N\Ff Q$ then we say that $Q$ {\it is normal in} $\Ff$. The largest subgroup of $S$ which is normal in $\Ff$ will be denoted $O_p(\Ff)$.

The {\it centralizer of $Q$ in $\Ff$} is the category $\C\Ff Q$ on $\C SQ$
having as morphisms those morphisms $\varphi\in\Hom_\Ff(R,T)$, with $R$ and $T$ subgroups of $C_S(Q)$, satisfying that there exists a morphism $\widehat{\varphi}\in\Hom_\Ff(QR,QT)$ such that $\widehat{\varphi}|_Q=\id_Q$ and $\widehat{\varphi}|_R=\varphi$.

Also denote by $\N SQ\C\Ff Q$ the category on $\N SQ$ having as morphisms all group homomorphisms $\varphi: P \rightarrow R$, for $P$ and $R$ subgroups of $\N SQ$, for which there exists a morphism $\psi: QP \rightarrow QR$ and $x \in\N SQ$ such that $\psi|_Q = c_x$ (the morphism induced by conjugation by $x$) and $\psi|_P=\varphi$.

\begin{prop}[\cite{puig:frobenius}, Proposition 2.8]\label{prpnorm}
If $Q$ is fully $\Ff$-normalized then $\N\Ff Q$ is a fusion system on $\N SQ$. 
If $Q$ is fully $\Ff$-centralized then $\C\Ff Q$ and $\N SQ\C \Ff Q$ are fusion systems on $\C SQ$.
\end{prop}

Alperin's theorem on $p$-local control of fusion also holds for fusion systems. First we set up this theorem's notations and terminology. If $\varphi\in\Aut_\Ff(S)$ we say that $\varphi$ is a maximal $\Ff$-automorphism. If $\varphi\in\Aut_\Ff(E)$, with $E$ an $\Ff$-essential subgroup of $S$, we say that $\varphi$ is an essential $\Ff$-automorphism. Alperin's fusion theorem asserts that the essential and maximal $\Ff$-automorphisms suffice to determine the whole fusion system $\Ff$.

\begin{thm}[Alperin]\label{AlperinsThm} Any morphism $\varphi\in\Hom_\Ff(Q,S)$ can be written as the composition of restrictions of essential $\Ff$-automorphisms, followed by the restriction of a maximal $\Ff$-automorphism.
More precisely, there exists
\begin{itemize}
\item[(a)] an integer $n\ge 0$,
\item[(b)] a set of $\Ff$-isomorphic subgroups of $S$, $Q=Q_0,\,Q_1,\,\dots,\,Q_n,\,Q_{n+1}=\varphi(Q)$,
\item[(c)] a set of $\Ff$-essential, fully $\Ff$-normalized subgroups
$E_i$ of $S$ containing $Q_{i-1}$ and $Q_i$, for all $1\le i\le n$,
\item[(d)] a set of essential automorphisms $\psi_i\in\Aut_\Ff(E_i)$ satisfying $\psi_i(Q_{i-1})=Q_i$, for all $1\le i\le n$ and
\item[(e)] a maximal automorphism $\psi_{n+1}\in\Aut_\Ff(S)$ satisfying $\psi_i(Q_n)=Q_{n+1}$,
\end{itemize}
such that we have
$$\varphi(u)=\psi_{n+1}\psi_{n}\ldots\psi_1(u),\text{ for all }u\in Q\,.$$
\end{thm}

The reader can find a proof of this theorem in \cite{stancu}; an alternative proof of this theorem in a different axiomatic setting was given by Puig \cite[Corollary 3.9]{puig:frobenius} and another in a less general form, using $\Ff$-centric, $\Ff$-radical subgroups instead of $\Ff$-essential subgroups, can be found in \cite[Theorem A.10]{blo}. We use this later form in the proof of Lemma \ref{Op_cr}\,. 

The classical examples of a fusion systems are the ones coming from the $p$-local structure of a finite group $G$. If $S$ is a Sylow $p$-subgroup of $G$ then we denote by $\Ff_S(G)$ the fusion system on $S$ having as morphisms
$$\Hom_{\Ff_S(G)}(P,Q)=\N G{P,Q}/\C GP =\Hom_G(P,Q)$$
where $P,\,Q$ are subgroups of $S$ and $\N G{P,Q}=\{g\in G\,|\,^gP\le Q\}$ is the $G$-transporter from $P$ to $Q$.

There are examples of fusion systems that do not come from a finite group (see eg. \cite{ruizviruel}). But there are also particular cases when one can construct a finite group with $p$-local structure equivalent to a given fusion system. This is the case for constrained fusion systems. The fusion system $\Ff$ is said to be {\it constrained} if $O_p(\Ff)$ is $\Ff$-centric. Any constrained fusion system was proven to come from a finite group by Broto, Castellana, Grodal, Levi and Oliver:

\begin{thm}\cite[Theorem 4.3]{bcglo1}\label{thmconstraint}
Let $\Ff$ be a fusion system on $S$ and suppose that there exists an $\Ff$-centric subgroup $Q$ of $S$
such that $\N\Ff Q=\Ff$ (in particular $\Ff$ is constrained). Then there exists a, unique up to isomorphism,
finite $p'$-reduced $p$-constrained group $L_Q^{\Ff}$ (i.e $O_{p'}(L_Q^{\Ff})=1$ and $Q$ is a normal 
subgroup of $L_Q^{\Ff}$) having $S$ as a Sylow $p$-subgroup and such that
$\Ff= \Ff_S (L_Q^{\Ff})$. Furthermore $L_Q^{\Ff}/Z(Q)\simeq\Aut_{\Ff}(Q)$.
\end{thm}

%%%%%%%%%%%%%%%%%%%%%%%%%%%%%%%%%%%%%%%%%%%%%%
\section{Further results on fusion systems} %%
%%%%%%%%%%%%%%%%%%%%%%%%%%%%%%%%%%%%%%%%%%%%%%

Let $G$ be a finite group and $p$ a prime divisor of its order. If $A \trianglelefteq B \leq G$ then $B/A$ is called a {\it section} of $G$. We say that $H$ {\it is involved in} $G$ if $H$ is isomorphic to a section of $G$. If $H$ is not involved in $G$ then $G$ is called $H$-{\it free}.

Following \cite{kl05} we say that the {\it fusion system $\Ff$ on $S$ is $H$-free} if $H$ is not involved in any of the groups $L_Q^{\Ff}$, for $Q$ running over the set of 
$\Ff$-centric, $\Ff$-radical and fully $\Ff$-normalized subgroups of $S$. In some particular cases the property of being $H$-free passes to subsystems and quotient systems as the next two results from \cite{kl05} show.

\begin{prop}\cite[Proposition 6.3]{kl05}\label{subsyst} Let $\Ff$ be a fusion system on a finite $p$-group $S$ and let $Q$ be a fully $\Ff$-normalized subgroup of $S$. If $\Ff$ is $H$-free, then so is any fusion subsystem of $\Ff$ which lies between $\N\Ff Q$ and $\N SQ\C\Ff Q$. In particular, if $\Ff$ is $H$-free, so are $\N \Ff Q$ and $\N SQ \C\Ff Q$.
\end{prop}

Let $\Ff$ be a fusion system on $S$ and let $Q$ be a subgroup of $S$ with the property $\Ff = \N \Ff Q$. The category $\Ff /Q$ on $S/Q$ is defined as follows: for $Q \leq P,R \leq S$, a group homomorphism $\psi: P/Q \rightarrow R/Q$ is a morphism in $\Ff/Q$ if there is a morphism $\varphi \in \Hom_\Ff(P,R)$ satisfying $\psi(xQ) = \varphi(x)Q$ for all $x \in P$. The fact that $\Ff/Q$ is a fusion system on $S/Q$ is due to Puig \cite{puig:notes}, see also \cite[Proposition 2.8]{kl05}.

\begin{prop}\cite[Proposition 6.4]{kl05}\label{quotient} Let $\Ff$ be a fusion system on a finite $p$-group $S$ and let $Q$ be a normal subgroup of $S$ such that $\Ff=\N \Ff Q$. If $\Ff$ is $H$-free then $\Ff /Q$ is also $H$-free.
\end{prop}

The following result generalizes the main technical step in the proof of \cite[Proposition 5.2]{kl05}.

\begin{prop}\label{wellplaced} Let $\Ff$ be a fusion system on a finite $p$-group $S$ and let $W_i$,
$1 \leq i \leq n$ be subgroups of $S$ such that
\begin{itemize}
\item[(a)] the subgroup $W_{i+1}$ is a characteristic subgroup of $\N S{W_i}$ for $1 \leq i \leq n-1$;
\item[(b)] the subgroup $W_i$ is fully $\Ff$-normalized for $1 \leq i \leq n-1$.
\end{itemize}
Then there exists a morphism $\varphi \in \Hom _\Ff (\N S{W_n}, S)$ such that $\varphi(W_i)$ is fully $\Ff$-normalized for all $1 \leq i \leq n$. In particular $\varphi(\N S{W_i})=\N S{\varphi (W_i)}$. If, moreover, $W_i$ is $\Ff$-centric and/or $\Ff$-radical for some $1 \le i \le n$, then so is $\varphi(W_i)$.
\end{prop}

\begin{proof}
If follows from Lemma $2.2(a)$ that there exists a morphism $\varphi\in\Hom _{\Ff} (N_S(W_n), S)$ such that $\varphi(W_n)$ is fully $\Ff$-normalized. According to condition (a), $W_{i+1}$ is a characteristic subgroup of $N_S(W_i)$, for $1 \leq i \le n-1$ and therefore $N_S(W_i) \le N_S(W_{i+1})$. In particular $N_S(W_i) \leq N_S(W_n)$ and the morphism $\varphi$ is defined on $\N S{W_i}$ for all $1 \leq i \leq n$.

Next we show $\varphi(W_i)$ is fully $\Ff$-normalized for all $1 \le i \le n-1$. By elementary group theory $\varphi(N_S(W_i)) \le N_S(\varphi(W_i))$. Since $\varphi$ is injective it follows that $|N_S(W_i)| \le |\varphi(N_S(W_i))|$. But, according to (b), $W_i$ is fully $\Ff$-normalized and $|N_S(W_i)| \ge |N_S(\varphi(W_i))|$. It follows now that $|N_S(W_i)| = |N_S(\varphi(W_i))|$ which shows that $\varphi(W_i)$ is fully $\Ff$-normalized and that $\varphi (N_S(W_i)) = N_S(\varphi(W_i))$, for all $1 \le i \le n-1$.

Since $W_i$ is fully $\Ff$-normalized for all $1 \leq i \leq n$, it is also fully $\Ff$-centralized. Thus $\varphi (\C S{W_i})=\C S{\varphi (W_i)}$ and if $W_i$ is $\Ff$-centric, then so is $\varphi(W_i)$. Moreover it is a general fact that $\,^\varphi\Aut_\Ff(W_i)=\Aut_\Ff(\varphi(W_i))$
and $\,^\varphi\Aut_S(W_i)=\Aut_S(\varphi(W_i))$ so if $W_i$ is $\Ff$-radical, then so is $\varphi(W_i)$.
\end{proof}

We need the following definition reproduced from \cite[Definition 5.1]{kl05}.

\begin{dfn}\label{charfunct}
A {\it positive characteristic functor} is a map sending any nontrivial finite $p$-group $S$ to a nontrivial characteristic subgroup $W(S)$ of $S$ such that $W(\varphi(S)) = \varphi(W(S))$ for every $\varphi \in \Aut(S)$. A positive characteristic functor is a {\it Glauberman functor} if whenever $S$ is a Sylow $p$-subgroup of a $Qd(p)$-free finite group $L$ which satisfies $C_L(O_p(L)) = Z(O_p(L))$, then $W(S)$ is normal in $L$.
\end{dfn}

Using Proposition \ref{wellplaced} we can give a different proof for Proposition $5.3$ in \cite{kl05}.

\begin{prop}\cite[Proposition 5.3]{kl05}\label{normal} Let $\Ff$ be a fusion system on a finite $p$-group $S$ and let $W$ be a positive characteristic functor. Assume that for any non-trivial, proper, $\Ff$-centric, $\Ff$-radical, fully $\Ff$-normalized subgroup $Q$ of $S$ the following holds $\N \Ff Q = \N {\N \Ff Q} {W(\N SQ)}$. Then $\Ff = \N \Ff {W(S)}$.
\end{prop}

\begin{proof}
Suppose that the conclusion does not hold. By Alperin's Fusion Theorem there exists a proper, $\Ff$-centric, $\Ff$-radical, fully $\Ff$-normalized subgroup $Q$ of $S$ such that $\Aut_{\N \Ff {W(S)}}(Q) \subsetneq \Aut_\Ff(Q)$.

Set $W_1=Q$ and define recursively $W_{i+1}:=W(\N S{W_i})$. So $W_{i+1}$ is characteristic in $\N S{W_i}$ implying that we get the following inclusions $\N S{W_i}< \N S{\N S {W_i}} \leq \N S{W_{i+1}}$ this in its turn implies there exists $n\ge 1$ such that the sequence of $\N S{W_i}$ for $1\le i\le n$ is strictly increasing and $\N S{W_n}=S$. Observe that if $\N S{W_i} = \N S{W_{i+1}}$ then $\N S{W_i} = \N S{\N S {W_i}} \leq S$ and by an elementary property of $p$-groups it follows that $N_S(W_i)=S$, thus indeed the sequence eventually reaches $S$.

Moreover the sequence $\{W_i,\,1\le i\le n+1\}$ can be chosen so that all its terms are fully $\Ff$-normalized. This can be done recursively by applying Proposition \ref{wellplaced}, for all $2\le k\le n+1$ to the partial subsequences $\{W_i,\,1\le i\le k\}$. The $W_i$'s are successively modified by replacing them with their images through the morphism $\varphi$ given by Proposition \ref{wellplaced}.

Consider the sequence of the normalizers in $\Ff$ of the $W_i$'s for $1\le i\le n+1$. Given that $W_i$, $1\le i\le n$ are fully $\Ff$-normalized, we have that $\N\Ff{W_i}$ is a fusion system on $\N S{W_i}$. It follows from our assumption that $\N\Ff{W_i}\subseteq\N\Ff{W_{i+1}}$ for all $1 \le i \le n-1$. But then $\N\Ff{Q}=\N\Ff{W_1}\subseteq\N\Ff{W_{n+1}}=\N\Ff{W(S)}$. At the level of morphisms on $Q$ this gives $\Aut_\Ff(Q)\subseteq\Aut_{\N\Ff{W(S)}}(Q)$ which is a contradiction with the initial supposition on $Q$.
\end{proof}

Proposition \ref{normal} is used in \cite{kl05} to prove the next important result.

\begin{prop}\cite[Proposition 3.4]{kl05}\label{quot} Let $S$ be a finite $p$-group and let $Q$ be a normal subgroup of $S$. Let $\Ff$ and $\Gg$ be fusion systems on $S$ such that $\Ff=S \C \Ff Q$ and such that $\Gg \subseteq \Ff$. Let $P$ be a normal subgroup of $S$ containing $Q$. We have $\Gg = \N \Ff P$ if and only if $\Gg /Q = \N {\Ff /Q} {P / Q}$.
\end{prop}

Next we give an application of the Frattini argument to fusion systems. The group theoretic result states that, if $G$ is a finite group, then thefactorization $G=\C GQ\N GR$ holds with $Q=O_p(G)$ and $R=\C G{Q\C SQ}$. This is easily seen to be true as $\C GQ \trianglelefteq G$ and $\N G{\C SQ} \leq \N GR$, then an application of the Frattini argument gives the result.

\begin{lem}\label{Op_cr} Let $\Ff$ be a fusion system on $S$, $Q=O_p(\Ff)$ and $R=Q\C SQ$. Set $\Ff_1=S\C \Ff Q$ and $\Ff_2=\N \Ff R$. Then $\Ff=\langle \Ff_1, \Ff_2 \rangle$.
\end{lem}

\begin{proof}
First remark that $\Ff$, $\Ff_1$ and $\Ff_2$ are fusion systems on $S$ with $\Ff$ containing the other two. By Alperin's fusion theorem (see Theorem \ref{AlperinsThm}), it is enough to prove that for every $\Ff$-centric, $\Ff$-radical, fully $\Ff$-normalized subgroup $U$ of $S$ we have $\Aut_\Ff(U)=\Aut_{\Gg}(U)$ with $\Gg= \langle \Ff_1, \Ff_2 \rangle$.

We shall prove that every $\varphi\in\Aut_\Ff(U)$ can be written as composition of morphisms in $\Ff_1$ and $\Ff_2$ and thus will be contained in $\Aut_{\Gg}(U)$. This will finish our proof as the opposite inclusion is clearly satisfied.

Given that $U$ is $\Ff$-centric and $\Ff$-radical we have by \cite[Proposition 5.6]{intro} that $Q\le U$. Hence $\varphi$ restricts to an automorphism $\theta\in\Aut_\Ff(Q)$. Now we have that $N_\theta$ contains $U$ and $R$ so it contains $UR$. Given that $Q$ is fully $\Ff$-normalized $\theta$ extends to $\chi\in\Hom_\Ff(UR,S)$. Moreover $\chi(R)=\chi(Q)\C S{\chi(Q)}=R$ so in fact $\chi\in\Hom_{\Ff_2}(UR,S)$.

Denote by $\psi$ the restriction to $U$ of $\chi$; then $\psi\in\Hom_{\Ff_2}(U,\psi(U))$.
Both $\varphi$ and $\psi$ restrict as $\theta$ on $Q$ so $\varphi\circ\psi^{-1}$ belongs to $\Hom_{\Ff_1}(\psi(U),U)$. The conclusion in the lemma follows as $\varphi=\varphi\circ\psi^{-1}\circ\psi\in\Aut_{\Gg}(Q)$.
\end{proof}

We close this section with a straightforward result on fusion control in fusion systems.

\begin{lem}\label{normalizerW} Let $W$ be a fully $\Ff$-normalized subgroup of $S$ and suppose that there are two fusion subsystems $\Ff_1$ and $\Ff_2$ of $\Ff$ such that $\Ff=\langle \Ff_1,\Ff_2 \rangle$. If moreover $\Ff_1 = \N {\Ff_1} W$ and $\Ff_2 = \N {\Ff_2} W$. Then $\Ff = \N \Ff W$.
\end{lem}

\begin{proof} We have that $\langle \N {\Ff_1}W, \N {\Ff_2}W \rangle \subseteq \N \Ff W \subseteq \Ff$. The result follows. \end{proof}

%%%%%%%%%%%%%%%%%%%%%%%%%%%%%%%%%%%%%%%%%%%%%%%%%%%%%%%%%%
\section{A characteristic subgroup of $S$}\label{charsg}%%
%%%%%%%%%%%%%%%%%%%%%%%%%%%%%%%%%%%%%%%%%%%%%%%%%%%%%%%%%%

Let $S$ be a finite $p$-group. In this section we construct a subgroup $W(S)$ which is characteristic in $S$ and such that $\Omega(Z(S)) \leq W(S) \leq \Omega(Z(J(S)))$, and with the property that $W(S) \trianglelefteq \Ee$ for all $(\varphi, \Ee) \in \Uu_J$, with $\Uu_J$ a class of embeddings defined below. The notation $\Omega (H)$ stands for the subgroup of $H$ generated by all the elements of order $p$, while $J(S)$ denotes the Thompson subgroup of $S$ defined in Section 1. We shall give below two different, although equivalent, constructions of this characteristic subgroup of $S$ which we shall denote $W(S)$ and $W$. The first construction follows the approach developed by Stellmacher \cite{st92,st96} for finite groups, in which such a subgroup is approximated from various subgroups of $Z(J(S))$. The second construction uses basic properties of fusion systems.

\subsection{The group W(S)} An {\it embedding} is a pair $(\varphi, \Ee)$ where $\varphi \in \Aut(S)$ and $\Ee$ is a category on $\varphi(S)= S$. Let $\Cc$ denote the family of all embeddings of $S$. A nonempty subclass $\Uu$ of $\Cc$ is {\it characteristically closed} if $(\varphi \alpha, \Ee) \in \Uu$ whenever $(\varphi, \Ee) \in \Uu$ and $\alpha \in \Aut(S)$.

An {\it equivalence} between two embeddings $(\varphi _1, \Ee_1)$ and $(\varphi _2, \Ee_2)$ is a morphism
$(\alpha, \Theta):\Ee_1 \rightarrow \Ee_2$ with $\alpha \varphi _1 = \varphi _2$ and
$\Hom_{\Ee_2}(\varphi _2(Q), \varphi _2(R)) = \alpha \circ \Hom _{\Ee_1}(\varphi _1(Q), \varphi_1(R)) \circ {\alpha ^{-1}}_{|\varphi_1 (Q)}$.
The equivalence of embeddings defines an equivalence relation on $\Uu$. Since $S$ is a finite group,
the collection of equivalence classes $[\Uu]$ is a finite set.

Let $O_S(\Uu)$ denote the largest subgroup of $S$ which satisfies the property that $\varphi (O_S(\Uu))$ is normal in $\Ee$ for every embedding $(\varphi, \Ee)$ in $\Uu$.

\begin{lem}Let $\Uu$ be a characteristically closed subclass of $\Cc$ and let $\alpha \in \Aut(S)$. Let $Q$ be a subgroup of $S$ with the property that $\varphi (Q) \trianglelefteq \Ee$ for every $(\varphi, \Ee) \in \Uu$. Then $\varphi(\alpha (Q)) \trianglelefteq \Ee$ for every $(\varphi, \Ee)$. In particular $O_S(\Uu)$ is a characteristic subgroup of $S$.
\end{lem}
\begin{proof} Observe that since $\Uu$ is characteristically closed and $(\varphi, \Ee) \in \Uu$ then $(\varphi\alpha, \Ee) \in \Uu$ and thus $\varphi (\alpha (Q)) \trianglelefteq \Ee$ for every $(\varphi, \Ee) \in \Uu$. The fact that $O_S(\Uu)$ is $\Aut(S)$-invariant follows from its definition.
\end{proof}

Let $\Uu_J$ denote the class of embeddings $(\varphi, \Ee)$ which satisfy the following conditions:
\begin{itemize}
\item[$(U1)$]$\Uu_J$ is characteristically closed.
\item[$(U2)$]$J(\varphi(S))=J(S)$ is normal in $\Ee$ for all $(\varphi, \Ee) \in \Uu_J$.
\item[$(U3)$]$\Ee$ is a $Qd(p)$-free fusion system.
\end{itemize}

For a $p$-group $P$ we set $A(P)=\Omega(Z(P))$ and $B(P)=\Omega(Z(J(P)))$. Remark that $A(P)\le B(P)$ as $Z(P)\le J(P)$. Note that $\alpha(A(P))=A(\alpha(P))$ and $\alpha(B(P))=B(\alpha(P))$ for all $\alpha \in \Aut(P)$, as $A(P)$ and $B(P)$ are characteristic subgroups of $P$.

Define recursively a subgroup $W(S) \leq B(S)$ as follows. Let
$$W_0:= A(S)=\Omega(Z(S)) \leq B(S)$$
and assume that for $i \geq 1$ the subgroups $W_0, W_1, \ldots W_{i-1}$ with
$$W_0 < W_1 <\ldots W_{i-1} \leq B(S)$$
are defined. If $\varphi(W_{i-1})\trianglelefteq \Ee$ for all $(\varphi, \Ee) \in \Uu_J$ then set $W(S):=W_{i-1}$. Otherwise, choose $(\varphi _i, \Ee_i) \in \Uu_J$ to be such that $\varphi_i(W_{i-1})$ is not normal in $\Ee_i$ and define
$$W_i:= \varphi^{-1}_i \left< \varphi_i (W_{i-1})^{\Ee_i} \right> = \varphi^{-1}_i \left< \psi(\varphi _i(W_{i-1})) | \psi \in \Hom_{\Ee_i}(\varphi_i(W_{i-1}),\varphi_i(S)) \right>$$
to be the preimage in $S$ of the group generated by the $\Ee_i$-orbit of $\varphi_i(W_{i-1})$.

Since $B(\varphi_i(S))$ is a characteristic subgroup of $J(\varphi_i(S))$, which is itself normal in $\Ee_i$, it follows that $B(\varphi_i(S))$ is also normal in $\Ee_i$. Clearly $\varphi_i(W_i) \leq B(\varphi_i(S))$ since $\varphi_i(W_i)$ is generated by various conjugates of $\varphi_i(W_{i-1})$ and $\varphi_i(W_{i-1}) \leq B(\varphi_i(S))$ by construction. Thus we have:
$$A(\varphi_i(S)) \leq \varphi_i(W_{i-1}) < \varphi_i(W_i) \leq B(\varphi_i(S)) \trianglelefteq \Ee_i$$ as $A(\varphi_i(S))=\varphi_i(A(S))$ and $B(\varphi_i(S))=\varphi_i(B(S))$ for $\varphi_i \in \Aut (S)$. Then it follows:
$$A(S) \leq W_{i-1} < W_i \leq B(S)$$
As $S$ is finite, this recursive definition terminates after a finite number $n$ of steps and $W(S) := W_n$. Therefore we obtain a chain of subgroups of $Z(J(S))$:
$$A(S) = W_0 < W_1 \ldots < W_i < \ldots < W_n=W(S) \leq B(S)$$
and $\varphi(W(S)) \trianglelefteq \Ee$ for all $(\varphi, \Ee) \in \Uu_J$.

The group $W(S)$ depends on $S$ only and it is independent of the pairs $(\varphi_i, \Ee_i)$. To see this assume that we defined in an analogous way
$$W_0 = \overline{W}_0<\overline{W}_1 < \ldots < \overline{W}_{\tilde{n}} =: \overline{W}(S)$$
for suitable pairs $(\overline{\varphi}_j, \overline{\Ee}_j)$ in $\Uu_J$ and for $j =1, \ldots \overline{n}$. First note that $\overline{W}_0=W_0 \leq W(S) \cap \overline{W}(S)$. Thus $\overline{\varphi}_1(\overline{W}_0)=\overline{\varphi}_1(W_0) \leq \overline{\varphi}_1(W(S)) = W(S)$ as $W(S) \trianglelefteq \overline{\Ee}_1$, and therefore $\overline{\varphi}_1 (\overline{W}_1) = \langle \overline{\varphi}_1 (\overline{W}_0) ^{\overline{\Ee}_1} \rangle \leq W(S)$ which implies $\overline{W}_1 \leq W(S)$. Proceed by induction on $j$; a similar argument shows that since $\overline{W}_{j-1} \leq W(S)$ then $\overline{\varphi}_j(\overline{W}_j) \leq W(S)$ and $\overline{W}_j \leq W(S)$. Therefore $\overline{W}(S) \leq W(S)$. Similarly $W(S) \leq \overline{W}(S)$ and thus $W(S) = \overline{W}(S)$.

\begin{lem}Let $\alpha \in \Aut(S)$. Then $W(\alpha(S)) = \alpha(W(S))$. In particular, $W(S)$ is a characteristic subgroup of $S$, nontrivial if $S$ is nontrivial.
\end{lem}

\begin{proof} The mapping $(\varphi, \Ee) \rightarrow (\varphi \alpha, \Ee)$ defines a bijection on $\Uu_J$. Under this map, the chain of subgroups:
$$A(S) = W_0 < \ldots < W_i < \ldots < W_n = W(S) \leq B(S)$$
is taken to the following chain:
$$A(S) = \alpha(W_0) < \alpha(W_1) < \ldots < \alpha(W_i) < \ldots <\alpha(W_n) = \alpha(W(S)) \leq B(S).$$
Therefore $W(\alpha(S))=\alpha(W(S))$. The last statement follows from the fact that $\Omega(Z(S))\leq W(S)$ and $Z(S) \ne 1$ if $S \ne 1$.
\end{proof}

\subsection{The group W} Denote by $\Cc_J$ the class of categories $\Ff$ on $S$ which satisfy the following conditions:
\begin{itemize}
\item[$(C1)$]$J(S)$ is normal in $\Ff$ for all $\Ff \in \Cc_J$.
\item[$(C2)$]$\Ff$ is a $Qd(p)$-free fusion system.
\end{itemize}

\begin{prop}Let $W_0=\Omega (Z(S))$ and define
$$W:=\langle \psi(W_0) \; | \; \psi \in \Hom_\Ff(J(S), S) \;{\rm for}\;\Ff \in \Cc_J \rangle.$$
The subgroup $W$ is a nontrivial characteristic subgroup of $S$.
\end{prop}

\begin{proof}For all $\alpha \in \Aut(S)$, we will show that $\alpha(W) = W$. Let $\Ff$ be a category on $S$.
Denote by $\Ff^\alpha$  the category on $S$ having as sets of morphisms
$$\Hom_{\Ff^\alpha}(Q,R)=\alpha^{-1}\circ\Hom_\Ff(\alpha(Q),\alpha(R))\circ\alpha\,.$$
Note that if $\Ff\in\Cc_J$ then $\Ff^\alpha\in\Cc_J$, and if $\psi\in\Hom_\Ff(Q,R)$ then $\alpha \psi \alpha^{-1} \in \Hom_{\Ff^{\alpha^{-1}}}(\alpha(Q),\alpha(R))$.
Thus we have:
$$\begin{array}{lll}
\alpha(W)&:=\langle \alpha\psi(W_0) \; | \; \psi \in \Hom_\Ff(J(S),S) \;{\rm for}\;\Ff \in \Cc_J \rangle = &\\
&= \langle \alpha\psi\alpha^{-1}(\alpha(W_0)) \; | \; \alpha\psi\alpha^{-1} \in \Ff^{\alpha^{-1}}(\alpha(J(S)),\alpha(S)) \;{\rm for}\;\Ff \in \Cc_J \rangle = &\\
&= \langle \widetilde\psi(W_0) \; | \; \widetilde\psi \in \Ff^{\alpha^{-1}}(J(S),S) \;{\rm for}\;\Ff \in \Cc_J \rangle &\\
&\leq W
\end{array}
$$
where in the last equality we use that $W_0$ and $J(S)$ are characteristic subgroups of $S$.
But since $\alpha$ is injective, it follows that $|W| =|\alpha(W)|$ and therefore $\alpha(W)=W$, proving that $W$ is characteristic in $S$.
\end{proof}

\subsection{Stellmacher functor} Given that $(\varphi, \Ff)$ and $(\alpha\varphi, \Ff^{\alpha^{-1}})$ are equivalent as embeddings, $\Ff\in\Cc_J$ if and only if $\Ff^\alpha\in\Cc_J$ for any $\varphi,\alpha\in\Aut(S)$, and $\Ff$ a category on $S$, the two definitions $W(S)$ and $W$ represent the same subgroup of $S$, that is $W=W(S)$. It follows from Lemma $4.2$ that the functor $S \rightarrow W(S)$, for $S$ a finite $p$-group, is a positive characteristic functor in the sense of Definition \ref{charfunct}. We shall call the functor $S \rightarrow W(S)$, with $W(S)$ constructed via one of the methods from this Section, a {\it Stellmacher functor}.

The Thompson subgroup of $S$, $J(S)$ is a characteristic, centric subgroup. Thus using $(C1)$, any $\Ff\in\Cc_J$ is a constrained fusion system on $S$ and by Theorem \ref{thmconstraint}, there exists a $p$-constrained finite group $L$ with $\Ff = \Ff_L(S)$ and satisfying the following conditions: $S$ is a Sylow $p$-subgroup of $L$; $\C L{O_p(L)} \leq O_p(L)$ and $L$ is $Qd(p)$-free. It follows from our construction that $W(S)$ is a characteristic subgroup of $S$ which is also normal in $L$. The construction of $W(S)$ depends on $S$ only, and the subgroup $W(S)$ is constructed in the same way  as Stellmacher does in the context of finite groups so it is the same subgroup of $S$. Finally, notice that $S \rightarrow W(S)$ is also a Glauberman functor, see Definition \ref{charfunct}.

%%%%%%%%%%%%%%%%%%%%%%%%%%%%%%%%%%%%%%%%%%%%%%%%%%%
\section{Proofs of the Theorems}%%%%%%%%%%%%%%%%%%%
%%%%%%%%%%%%%%%%%%%%%%%%%%%%%%%%%%%%%%%%%%%%%%%%%%%

\subsection*{Proof of Theorem 1.1}
Let $S$ be a finite $2$-group and let $\Ff$ be a $\Sigma_4$-free fusion system on $S$. Also $W(S)$ is the characteristic subgroup of $S$ defined in Section $\ref{charsg}$. We prove that the normalizer $\N\Ff{W(S)}$ of 
$W(S)$ in $\Ff$ is equal to $\Ff$.

This is true for the smallest fusion system on $S$, which is $\Ff_S(S)$. Suppose now by induction that all proper $\Sigma_4$-free subsystems and all $\Sigma_4$-free quotient systems $\Ff/Q$, with $Q$ a nontrivial normal subgroup of $\Ff$ satisfy Theorem \ref{mainS4}.

If $O_2(\Ff)=1$ then for every non-trivial fully $\Ff$-normalized subgroup $P$ of $S$ we have that $\N\Ff P$ is a proper subsystem of $\Ff$ (otherwise $P\le O_2(\Ff)$). But then $\N\Ff P$ satisfies Theorem \ref{mainS4} by induction as it is $\Sigma_4$-free by Proposition \ref{subsyst}. Hence $\N\Ff P=\N{\N\Ff P}{W(\N SP)}$ for every non-trivial fully $\Ff$-normalized subgroup $P$. An application of Proposition \ref{normal} gives now that $\Ff=\N\Ff{W(S)}$. So we can suppose $\bf(H1)$: $O_2(\Ff)\ne 1$.

Set $Q:=O_2(\Ff)$ and $R:=Q\C SQ$. These subgroups of $S$ are both nontrivial by $(H1)$.
If $Q=R$ then $Q$ is $\Ff$-centric. Consequently $\Ff$ is a constrained fusion system. According to \cite[Proposition 4.3]{bcglo} there exists a $2'$-reduced $2$-constrained finite group $L_Q$, which is an extension of $\Aut_\Ff(Q)$ by $Z(Q)$ and having $S$ as a Sylow $2$-subgroup. Thus $\Ff=\Ff_S(L_Q)$. Since $\Ff$ is $\Sigma_4$-free, the group $L_Q$ is $\Sigma_4$-free, by definition and given the $L_Q$ is $2$-constrained we have
$C_{L_Q}(O_2(L_Q))\le O_2(L_Q)$. Then, according to Stellmacher's main theorem in \cite{st96}, see also Section 1, the group $W(S)$ is normal in $L_Q$. This in its turn implies that $W(S)\trianglelefteq \Ff_S(L_Q)$ and therefore $\Ff=\N\Ff{W(S)}$. Thus we can also make the assumption $\bf(H2)$: $Q\ne R$
implying moreover that $\N\Ff R$ is a proper subsystem of $\Ff$.

Next we will see that we also have $\bf(H3)$: $S\C\Ff Q \ne \Ff$. Indeed suppose that $S\C\Ff Q=\Ff$. Then $\Ff/Q$ is a proper quotient system of $\Ff$ which is $\Sigma_4$-free by Proposition \ref{quotient}. The induction hypothesis gives now $\Ff/Q=\N{\Ff/Q}{W(S/Q)}$. Next, Proposition \ref{quot} gives that $\Ff=\N\Ff U$ where $U$ is the preimage in $S$ of $W(S/Q)$. As $U \trianglelefteq \Ff$ it follows that $U\le Q$, but this leads to a contradiction given that $W(S/Q)\ne 1$.

According to Lemma \ref{Op_cr}, we have $\Ff = \langle \Ff_1, \Ff_2 \rangle$ with $\Ff_1=S \C \Ff Q$ and $\Ff_2 = \N \Ff R$. By $(H2)$ and $(H3)$ both $\Ff_1$ and $\Ff_2$ are proper subsystems of $\Ff$, the induction hypothesis gives that $W(S) \trianglelefteq \Ff_1$ and that $W(S) \trianglelefteq \Ff_2$.

Notice that $W(S)$ is fully $\Ff$-normalized. By Lemma \ref{fullynormalized}, there exists a morphism $\varphi\in\Hom_\Ff(N_S(W(S)),S)$ with $\varphi(W(S))$ fully $\Ff$-normalized. As $N_S(W(S))=S$ and since $W(S)$ is characteristic in $S$ it follows that $W(S)=\varphi(W(S))$ is fully $\Ff$-normalized.

Finally, an application of Lemma \ref{normalizerW} gives the result: $\Ff=\N\Ff{W(S)}$.

\subsection*{Proof of Theorem 1.2}
Let $S$ be a finite $p$-group. Recall that the construction $W(S)$ described in Section $4$ and which associates to $S$ a nontrivial characteristic subgroup $W(S)$ gives rise to a Glauberman functor.

Assume now that $\Ff$ is a $Qd(p)$-free fusion system on $S$. If $p$ is an odd prime, it follows from \cite[Theorem B]{kl05} that $\Ff = \N\Ff{W(S)}$. If $p=2$ then $Qd(2)=\Sigma_4$ and the result is given by Theorem \ref{mainS4}.

\subsection*{Proof of Theorem 1.3}
Let $p$ be an odd prime. Let $\Ff$ be a fusion system over a finite $p$-group $S$. Let $W(S)$ be the characteristic subgroup of $S$ given by the Stellmacher functor. Since $\Ff_S(S)\subseteq \N \Ff {W(S)} \subseteq \Ff$ it is enough to show that if $\N \Ff{W(S)}=\Ff_S(S)$ then $\Ff = \Ff_S(S)$. The proof is similar to that of Theorem $A$ in \cite{kl05}; for the sake of completeness we will provide the details.

Let $\Ff$ be a minimal counterexample to Theorem $1.3$; thus $\N \Ff {W(S)}=\Ff_S(S)$ but $\Ff \ne \Ff_S(S)$, and all the proper subsystems and quotient systems of $\Ff$ satisfy Theorem $1.3$. Under this assumption we show that $\Ff$ is a constrained fusion system by proving that $Q :=O_p(\Ff)$ is a nontrivial $\Ff$-centric proper subgroup of $S$. This is attained in the following six steps.

{\it Step 1 : Any fusion system $\Gg$ on $S$ which is properly contained in $\Ff$ is equal to $\Ff_S(S)$}.

As $\Gg \subset \Ff$ it follows that $\N \Gg {W(S)} \subseteq \N \Ff {W(S)} = \Ff_S(S)$. But $W(S)$ is a characteristic subgroup of $S$ and therefore $\Ff_S(S) \subseteq \N \Gg {W(S)}$. Thus $\N \Gg {W(S)} = \Ff_S(S)$ and the minimality assumption on $\Ff$ implies that $\Gg = \Ff_S(S)$.

{\it Step 2 : Let $P$ be a fully $\Ff$-normalized subgroup of $S$ and set $A=N_S(P)$. Then there is $\varphi \in \Hom_\Ff(A,S)$ such that both $\varphi(P)$ and $\varphi(W(A))$ are fully $\Ff$-normalized}.

By Lemma $2.2(a)$, there is a morphism $\varphi:N_S(W(A)) \rightarrow S$ such that $\varphi(W(A))$ is fully $\Ff$-normalized. Since $W(A)$ is a characteristic subgroup of $A$, we have $N_S(P)=A \leq N_S(A) \leq N_S(W(A))$ and the morphism $\varphi$ can be restricted to $\varphi:A \rightarrow S$. According to Lemma $2.2(b)$, the group $\varphi(P)$ is also fully $\Ff$-normalized.

{\it Step 3 : The subgroup $Q = O_p(\Ff)$ is nontrivial}.

Recall that it is assumed that $\Ff_S(S) \subset \Ff$. Alperin's fusion theorem implies that there is a fully $\Ff$-normalized subgroup $P$ of $S$ with $\Ff_A(A) \subset \N \Ff P$, for $A = N_S(P)$. Choose the subgroup $P$ such that:\\
\hspace*{.4cm}a) $W(A)$ is fully $\Ff$-normalized;\\
\hspace*{.4cm}b) $N_S(P)=A$ has maximal order among subgroups $T$ with $\Ff_{N_S(T)}(N_S(T)) \subset \N \Ff T$.\\
The choice of $P$ and the fact that $A$ is a proper subgroup of $N_S(W(A))$ implies that $\N \Ff {W(A)} = \Ff_{N_S(W(A))}(N_S(W(A))$. Therefore $N_{\N \Ff P}(W(A))=\Ff_A(A)$.

If $\N \Ff P \subset \Ff$ then the minimality assumption on $\Ff$ implies that $\N \Ff P = \Ff _R(R)$, which contradicts our choice of $P$. Thus we have $\N \Ff P = \Ff$ and $1 \ne P \trianglelefteq \Ff$. Hence $1 \ne P \leq Q$ which proves that $Q \ne 1$.

{\it Step 4 : $Q$ is a proper subgroup of $S$}.

If $Q=S$ then $\Ff = \N \Ff S=\N \Ff {W(S)} = \Ff_S(S)$ contradicting our assumption on $\Ff$.

{\it Step 5 : $S \C\Ff Q = \Ff_S(S)$ when $Q = O_p(\Ff)$}.

We have $S \C\Ff Q \subseteq \Ff$. If $S \C\Ff Q \subset \Ff$ then {\it Step 1} implies that $S \C\Ff Q = \Ff_S(S)$ and we are done. Assume now that $S \C \Ff Q = \Ff$ and recall that $\Ff \ne \Ff_S(S)$. An application of Proposition \ref{quot}, with $\Gg=\Ff_S(S)$ and $\Ff=\N\Ff Q$, gives that $\Ff /Q \ne \Ff_S(S)/Q = \Ff_{S/Q} (S/Q)$. By {\it Step 3} the subgroup $Q$ is nontrivial and the minimality assumption on $\Ff$ implies that $\N {\Ff/Q} {W(S/Q)} \ne \Ff_{S/Q}(S/Q)$. Let $P$ be the inverse image of $W(S/Q)$ in $S$. Notice that $W(S/Q) \ne 1$, by the definition of $W(S)$, and thus $P$ properly contains $Q$. Also $P \trianglelefteq S$ and $\N {\Ff/Q} {W(S/Q)} = \N {\Ff/Q} {P/Q}$. Another application of Proposition \ref{quot} gives that $\N \Ff P \ne \Ff_S(S)$. Since $\N \Ff P \subseteq \Ff$, {\it Step 1} implies that $\N \Ff P = \Ff$ which is a contradiction with the fact that $P$ contains $Q$ properly.

{\it Step 6 : The subgroup $Q$ is $\Ff$-centric}.

If $Q = R=SC_S(Q)$ then $Q$ is $\Ff$-centric and we are done. So let us assume that $Q <R$. Notice that $R \trianglelefteq S$. Then $\N \Ff R$ is a proper subsystem of $\Ff$ and an application of {\it Step 1} gives that $\N \Ff R = \Ff_S (S)$. Recall also that by the previous step, $S \C \Ff Q = \Ff_S(S)$. Therefore, Lemma \ref{Op_cr} implies that $\Ff = \Ff_S(S)$, which is a contradiction to our choice of $\Ff$. Thus we must have $Q = R$.

Since $Q = O_p(\Ff)$ is a nontrivial normal centric subgroup of $\Ff$, the fusion system $\Ff$ is constrained. But this means by \cite[4.3]{bcglo1}, that there is a $p'$-reduced $p$-constrained finite group $L$ with $S$ as a Sylow $p$-subgroup and such that $Q=O_p(L)$. Furthermore $\Ff = \Ff_S(L)$ and therefore $\N \Ff {W(S)} = \Ff_S (\N L{W(S)}$.

Since $\N \Ff{W(S)}=\Ff_S(S)$ it follows that $\N L{W(S)}$ has a $p$-complement; see Remark $6.10$ in Appendix. According to the normal $p$-complement theorem of Thompson, $6.11$ below,  it follows that $L$ has a $p$-complement. Therefore $\Ff_S(S)=\Ff_S(L)$ and we reached a contradiction with our assumption on $\Ff$. This concludes the proof of the Theorem $1.3$.

%%%%%%%%%%%%%%%%%%%%
\section{Appendix}%%
%%%%%%%%%%%%%%%%%%%%

Let $p$ be an odd prime, $G$ a finite group and $S$ a Sylow $p$-subgroup of $G$. We say that $G$ is $p$-{\it stable} if and only if for every $p$-subgroup $Q$ of $G$ and every element $x$ of $\N GQ$ such that $[Q,x,x]=1$, we have that $x\C GQ \in O_p(\N GQ/\C GQ)$.

A classic result of special significance to the theory of finite groups is Glauberman's ZJ-theorem \cite{gl68}:

\begin{thma}[Glauberman] Let $p$ be an odd prime. Let $G$ be a finite, $p$-stable group such that $\C G{O_p(G)} \leq O_p(G)$. Then $Z(J(S))$ is a normal subgroup of $G$.
\end{thma}

Using the following:

\begin{propa}[14.7, \cite{gl71}] Assume that $p$ is odd and that $G$ is a finite group. Then the following conditions on $G$ are equivalent:\\
\hspace*{.5cm}(a) the group $Qd(p)$ is not involved in $G$;\\
\hspace*{.5cm}(b) every section of $G$ is $p$-stable.
\end{propa}

the $ZJ$-theorem can be reformulated as follows:

\begin{thma}[Glauberman] Let $p$ be an odd prime and let $G$ be a $Qd(p)$-free finite group with $\C G{O_p(G)} \leq O_p(G)$. Then $Z(J(S))$ is a normal subgroup of $G$.
\end{thma}

For $p=2$ the $ZJ$-theorem does not hold anymore; see \cite[Section 11]{gl71}. As noted by Glauberman \cite{gl71} a necessary and sufficient condition for every section of $G$ to be $2$-stable is that $G$ have a normal $2$-complement, which is too strong to be useful.

In a couple of papers \cite{st92,st96}, Stellmacher proved an analogous version of Glauberman's $ZJ$-theorem, by constructing a characteristic subgroup $W(S)$ of $S$ and extending the result for $p=2$. An overview of his method, including a sketch of the proof for the odd prime case can be found in \cite[Section 9.4]{ks04}. The main theorem in \cite{st96} reads as follows:

\begin{thma}[Stellmacher] Let $S$ be a nontrivial finite $2$-group. Suppose that $G$ is a finite group satisfying the following:\\
\hspace*{.4cm}  (I) $G$ is $\Sigma_4$-free,\\
\hspace*{.4cm} (II) $S \in \Syl_2(G)$ and $C_G(O_2(G)) \leq O_2(G)$,\\
\hspace*{.4cm}(III) Every non-abelian simple section of $G$ is isomorphic to $Sz(2^{2n+1})$ or $PSL_2(3^{2n+1})$.
Then there exists a nontrivial characteristic subgroup $W(S)$ of $S$ which is normal in $G$.
\end{thma}

Next, consider a couple of useful lemmas:

\begin{lema}[Chp. II, Lemma 2.3, \cite{gl77}] The following conditions are equivalent:\\
\hspace*{.5cm}(a) $\Sigma_4$ is involved in $G$;\\
\hspace*{.5cm}(b) There exists a $2$-subgroup $Q$ of $G$ such that $\Sigma_3$ is involved in $\N G Q/\C GQ$.
\end{lema}

\begin{lema}[Chp. II, Corollary 7.3, \cite{gl77}] Let $G$ be a non-abelian simple group. The following are equivalent:\\
\hspace*{.5cm}(a) $G$ is $\Sigma_3$-free;\\
\hspace*{.5cm}(b) $G$ is isomorphic to $Sz(2^{2n+1})$ or $PSL(2,3^{2n+1})$.
\end{lema}

\begin{rem} Note that if $G$ is $\Sigma_3$-free then $G$ is $\Sigma_4$-free. A finite group $G$ with $\C G{O_2(G)} \leq O_2(G)$ is $\Sigma_4$-free if and only if $G/O_2(G)$ is $\Sigma_3$-free \cite{st96}.
\end{rem}

Using the previous two lemmas and remark, we can rephrase Stellmacher's Theorem $6.4$ as follows:

\begin{thma}[Stellmacher] Let $S$ be a finite nontrivial $2$-group. Then there exists a nontrivial characteristic subgroup $W(S)$ of $S$ which is normal in $G$, for every finite $\Sigma_4$-free group $G$ with $S$ a Sylow $2$-subgroup and $\C G{O_2(G)} \leq O_2(G)$.
\end{thma}

If $G=S O_{p'}(G)$, with $S$ a Sylow $p$-subgroup of $G$, we say that $G$ {\it has a normal p-complement}. A standard result due to Frobenius (see \cite[8.6]{gl71} for example) is given below:

\begin{thma}[Frobenius] The following conditions are equivalent for a finite group $G$ with Sylow $p$-subgroup $S$:\\
\hspace*{.4cm}(a) $G$ has a normal $p$-complement;\\
\hspace*{.4cm}(b) if $Q$ is a non-identity subgroup of $G$ then $\N G Q /\C G Q$ is a $p$-group;\\
\hspace*{.4cm}(c) if $Q$ is a non-identity $p$-subgroup of $G$ then $\N G Q$ has a normal $p$-complement;\\
\hspace*{.4cm}(d) if two elements of $S$ are conjugate in $G$, they are conjugate in $S$.
\end{thma}

\begin{rem} The equivalence (a) $\Leftrightarrow$ (d) in the above theorem, states that $G$ has a normal $p$-complement if and only if $S$ controls fusion in $G$. In the language of fusion systems, $S$ controls $G$ fusion if and only if $\Ff_S(S)=\Ff_S(G)$.
\end{rem}

For odd primes Frobenius' result was improved by a result of Thompson. We give below a version of Thompson's p-complement theorem which uses Stellmacher's characteristic subgroup $W(S)$:

\begin{thma}\cite[9.4.7]{ks04} Let $G$ be a group, $p$ an odd prime, and $S \in \Syl_p(G)$. Then $G$ has a normal $p$-complement if and only if $\N G{W(S)}$ has a normal $p$-complement.
\end{thma}

\end{document}